\count100= 1
\count101= 10

\overfullrule 0pt

\magnification\magstep1
\hsize 4.43in
\vsize 7.28in



\font\bfab = cmbx9
\font\fa=cmr17
\font\fb=cmr12

\font\fab=cmr9
\font\sfab=cmr7
\font\fs=cmr6
\font\fd=cmr5
\font\slab=cmsl9
\font\tab= cmmi9
\font\sab= cmmi7
\font\ssab= cmmi6
\font\abst= cmsy9
\font\sabs= cmsy7
\font\ssabs= cmsy6
\font\itab= cmti9

\def\figfont{
    \textfont0 = \fab
    \scriptfont0 = \sfab
    \scriptscriptfont0 = \fs
    \textfont1 = \tab
    \scriptfont1 = \sab
    \scriptscriptfont1 = \ssab
    \textfont2 = \abst
    \scriptfont2 = \sabs
    \scriptscriptfont2 = \ssabs
    \let \sl = \slab
    \let \bf = \bfab
    \let \it = \itab
    \baselineskip 9pt
    \fab}


\def\CC{{\rm C\kern-.18cm\vrule width.6pt height 6pt depth-.2pt
\kern.18cm}}

\def\NN{{\mathop{{\rm I}\kern-.2em{\rm N}}\nolimits}}

\def\PP{{\mathop{{\rm I}\kern-.2em{\rm P}}\nolimits}}

\def\RR{{\mathop{{\rm I}\kern-.2em{\rm R}}\nolimits}}

\def\RRt{{\fa I}\kern-.2em{\fa R}}


\def\ZZ{{\mathop{{\rm Z}\kern-.28em{\rm Z}}\nolimits}}





\def\makebold#1{\mathord{\setbox0=\hbox{#1}%
       \copy0\kern-\wd0%
       \raise\dimen1\copy0\kern-\wd0%
       {\advance\dimen1 by \dimen1\raise\dimen1\copy0}\kern-\wd0%
       \kern\dimen0\raise\dimen1\copy0\kern-\wd0%
       {\advance\dimen1 by \dimen1\raise\dimen1\copy0}\kern-\wd0%
       \kern\dimen0\raise\dimen1\copy0\kern-\wd0%
       {\advance\dimen1 by \dimen1\raise\dimen1\copy0}\kern-\wd0%
       \kern\dimen0\raise\dimen1\copy0\kern-\wd0%
       \kern\dimen0\box0}}


\newbox\maboite

\def\boxit#1#2{\setbox\maboite=\hbox{\kern#1{#2}\kern#1}%
    \dimen1=\ht\maboite \advance\dimen1 by #1 \dimen2=\dp\maboite
\advance\dimen2 by #1%
    \setbox\maboite=\hbox{\vrule height\dimen1
depth\dimen2\box\maboite\vrule}%
    \setbox\maboite=\vbox{\hrule\box\maboite\hrule}%
    \advance\dimen1 by .4truept \ht\maboite=\dimen1%
    \advance\dimen2 by .4truept \dp\maboite=\dimen2 \box\maboite\relax}




%

\def\frac#1#2{{#1 \over #2}}


\def\ms{\medskip}

\def\noin{\noindent}










\def\pf{\noindent{\bf Proof: }}

\def\eop{\makeblanksquare6{.4}\ms}

\def\makeblanksquare#1#2{
\dimen0=#1pt\advance\dimen0 by -#2pt
      \vrule height#1pt width#2pt depth0pt\kern-#2pt
      \vrule height#1pt width#1pt depth-\dimen0 \kern-#1pt
      \vrule height#2pt width#1pt depth0pt \kern-#2pt
      \vrule height#1pt width#2pt depth0pt
}



\def\abstract#1{\bigskip\bigskip\medskip%
    {\narrower \baselineskip 9pt \fab \noindent {\bf Abstract.~~}%
    \textfont0 = \fab
    \scriptfont0 = \sfab
    \scriptscriptfont0 = \fs
    \textfont1 = \tab
    \scriptfont1 = \sab
    \scriptscriptfont1 = \ssab
    \textfont2 = \abst
    \scriptfont2 = \sabs
    \scriptscriptfont2 = \ssabs
    \let \it = \itab
    \let \sl = \slab
    #1\bigskip}\medskip}
\def\author#1{\bigskip\bigskip\centerline{\fb #1}}

\def\copyright{\hbox{{\fb o}\kern-.61em \raise .46ex \hbox{\fd c}}}
\def\footnoterule{\kern -3pt \hrule width 0truein \kern 2.6pt}
\def\leftheadline{\ifnum\pageno=\count100 \hfill%
  \else\rm\folio\hfil\it\shortauthor\fi}
\def\rightheadline{\ifnum\pageno=\count100 \hfill%
  \else\it\shorttitle\hfil\rm\folio\fi}
\def\title#1{\centerline {\fa #1}}
\def\titwo#1{\medskip \centerline {\fa #1}}

\def\titexp#1#2{\hbox{{\fa #1} \kern-.25em \raise .90ex \hbox{\fb #2}}\/}
\def\titsub#1#2{\hbox{{\fa #1} \kern-.25em \lower .60ex \hbox{\fb #2}}\/}

\nopagenumbers
\headline{\ifodd\pageno\rightheadline \else\leftheadline\fi}
\footline{\hfil}
\null\vskip 18pt
\centerline{}
\pageno=\count100
\count102=\count100
\advance\count102 by -1
\advance\count102 by \count101




\def\sect#1{\goodbreak\bigskip\smallskip\centerline{\bf\S #1}\medskip
    \noindent\ignorespaces}


\def\subsect#1{\goodbreak\bigskip\leftline{\bf#1}\medskip
  \noindent\ignorespaces}




\def\Address{\nonfrenchspacing\goodbreak\bigskip\obeylines}

\def\Remark#1.{\goodbreak\medskip\noin {\bf Remark#1.}}
\def\Example#1.{\goodbreak\medskip\noin {\bf Example#1.}}


\def\ref{\smallskip\global\advance\refnum by 1 \item{\the\refnum.}}
\newcount\refnum \refnum = 0

\def\References{\goodbreak\bigskip\centerline{\bf References}\bigskip
   \frenchspacing}


\def\ATA{Approx.\ Theory Appl.}

\def\JAT{J. Approx.\ Theory}



\title{Near-best univariate spline discrete quasi-interpolants}
\titwo{on non-uniform partitions}

\author{ D. Barrera, M.J. Iba\~nez, P. Sablonni\`ere, D. Sbibih}

\def\shorttitle{Univ Spline QI/non-uniform partitions}

\def\shortauthor{ps}


\abstract{Univariate spline discrete quasi-interpolants (abbr. dQIs)
are approximation operators using B-spline expansions with coefficients
which are
linear combinations of discrete values of the function to be approximated.
When working with nonuniform partitions, the main challenge is to find dQIs
which have
both good approximation orders and bounded uniform norms independent of the
given
partition. Near-best dQIs are obtained by minimizing an upper
bound of the infinite norm of dQIs depending on a certain number of free
parameters, thus reducing this norm. This paper is devoted to the study of
some families of near-best dQIs of approximation order 2.}


\sect{1.Introduction}

\noin
A spline quasi-interpolant (abbr. QI) of $f$ has the general form
$$
Qf=\sum_{\alpha\in A} \mu_{\alpha}(f) B_{\alpha}
$$
where $\{B_{\alpha},\alpha\in A\}$ is a family of B-splines forming a
partition of unity and
$\{\mu_{\alpha}(f),\alpha\in A\}$
is a family of linear functionals which are local in the sense that they
only use values of $f$
in some neighbourhood of  $\Sigma_{\alpha}=supp(B_{\alpha})$. The main
interest of QIs is that they
provide good approximants of functions without solving any linear system of
equations.
In the literature, one can find the three following types of QIs:
\ms
\noin
(i) {\sl Differential QIs} (abbr. DQIs) : the linear functionals are {\sl
linear combinations of
values of derivatives} of $f$ at some point in $\Sigma_{\alpha}$ (see e.g.
[5-7]).

\noin
(ii) {\sl Discrete QIs} (abbr. dQIs) : the linear functionals are {\sl
linear combinations of
values} of $f$ at some points in some neighbourhood of $\Sigma_{\alpha}$
(see e.g. [1-3], [6],
[9], [11], [13], [15-16], [24]).

\noin
(iii) {\sl Integral QIs} (abbr. iQIs) : the linear functionals are
{\sl linear combinations of weighted mean values} of $f$ in some
neighbourhood of
$\Sigma_{\alpha}$ (see e.g. [2-3], [6], [13-14], [24-25]).
\ms
\noin
In this paper and a subsequent one, we shall study various types of
univaraite dQIs and iQIs,
more specifically those that we call {\sl near-best} QIs which are defined
as follows:
\ms
\noin
(dQIs) assume that $\mu_{\alpha}(f)=\sum_{\beta\in F_{\alpha}}
\lambda_\alpha(\beta)
f(x_\beta)$ where the finite set of points $\{x_\beta, \beta\in
F_{\alpha}\}$ lies in some
neighbourhood of $\Sigma_{\alpha}$. Then it is clear that, for
$\Vert f \Vert_{\infty}\le 1$ and $\alpha\in A$,
$\vert \mu_{\alpha}(f) \vert \le \Vert \lambda_\alpha \Vert_1$,
where $\lambda_\alpha$ is the vector with components
$\lambda_\alpha(\beta)$, from which we deduce
immediately
$$
\Vert Q \Vert_{\infty} \le \sum_{\alpha\in A} \vert \mu_{\alpha}(f) \vert
B_{\alpha}\le
\max_{\alpha\in A}\vert \mu_{\alpha}(f) \vert\le \max_{\alpha\in A}\Vert
\lambda_\alpha
\Vert_1=\nu_1(Q).
$$
Now, assuming that $n=card(F_\alpha)$ for all $\alpha$, we can try to find a
$\lambda_\alpha^*\in\RR^n$ solution of  the minimization problem
$$
\Vert \lambda_\alpha^*\Vert_1=\min\{\Vert \lambda_\alpha\Vert_1;\;
\lambda_\alpha\in\RR^n, \;
V_\alpha \lambda_\alpha=b_\alpha\}
$$
where the linear constraints express that $Q$ is exact on some subspace of
polynomials. Thus, we finally
obtain
$$
\Vert Q \Vert_{\infty} \le \nu_1^*(Q)=\max_{\alpha\in A}\Vert
\lambda_\alpha \Vert_1.
$$

\noin
(iQIs) assume that $\mu_{\alpha}(f)=\sum_{\beta\in F_{\alpha}}
\lambda_\alpha(\beta)
\int_{\Sigma_\beta} M_{\beta}(t)f(t)dt$, where the B-splines $M_{\beta}$
are normalized by $\int
M_{\beta}=1$. Once again, for $\Vert f \Vert_{\infty}\le 1$, we have
$$
\vert \mu_{\alpha}(f)\vert\le \sum_{\beta\in F_{\alpha}} \vert\lambda_\alpha(\beta)\vert
\vert\int_{\Sigma_\beta} M_{\beta}(t)f(t)dt \vert \le \sum_{\beta\in F_{\alpha}}
\vert\lambda_\alpha(\beta)\vert=\Vert \lambda_\alpha \Vert_1
$$
whence, as we obtained above for dQIs,
$$
\Vert Q \Vert_{\infty} \le \max_{\alpha\in A} \Vert \lambda_\alpha
\Vert_1=\nu_1^*(Q).
$$
As emphasized by de Boor (see e.g. [5], chapter XII), a QI defined on non
uniform partitions has
to be {\sl uniformly bounded independently of the partition} (abbr. UB) in
order to be
interesting for applications. Therefore, the aim of this
paper is to define some families of dQIs satisfying this property and
having the {\sl smallest possible
norm}. As in general it is difficult to minimize the true norm of the
operator, {\sl we have
chosen to solve the minimization problems} defined above. A further paper
[25] will develop the
case of iQIs on nonuniform partitions. A few results are given in [24].
\ms
\noin
The paper extends some results of [1][13] and is organized as follows. We
first recall some
"classical" QIs of various types and we verify that they are UB. Then we
define and study several
families of discrete and integral QIs, depending on a finite number of
parameters, for which we
can find $\nu_1^*(Q)$. We show that this problem has always a solution (in
general non unique). We
give more specific examples for quadratic and cubic splines. Of particular
interest are the results
of theorems 3,5 and 6 where we show that some families of dQIs are
uniformly bounded independently of the
partition. Finally, we briefly give some applications to the
approximation of functions, to quadrature formulas and to pseudo-spectral
methods (see e.g. [12],
[29]). A parallel study of spline QIs is done in [2] for uniform partitions
of the real line and
in [3] for some uniform triangulations of the plane.

----------------

\sect{2. Notations}
\noin
We shall use classical B-splines of degree $m$ on a bounded interval
$I=[a,b]$ or on $I=\RR$. For the sake of simplicity, in the case $I=\RR$, we
take an increasing sequence of knots $T=\{t_i, i\in\ZZ\}$. In the  case
$I=[a,b]$, we take
the usual sequence $T$ of knots defined by (see [5][11][19][28]):
$$
a=t_{-m}=\ldots=t_0,\;\; b=t_n=\ldots=t_{n+m}
$$
$$
a<t_1<t_2<\ldots<t_{n-1}<b
$$
For $J=\{0,\ldots,n+m-1\}$,
the family of B-splines $\{B_j,j\in J\}$,
with support $\Sigma_j=[t_{j-m},t_{j+1}]$ is a basis of
the space
$S_m(I,T)$ of splines of degree $m$ on the interval $I$ endowed with the
partition $T$.
These B-splines form a partition of unity, i.e. $\sum_{j\in J} B_j=1$. We
denote $h_i=t_i-t_{i-1}$
for all indices $i$.

\noin
Let $\NN_m=\NN\cap [0,m-1]$ and $T_j=\{t_{j-r}, r\in \NN_m\}$:
we recall that the {\sl elementary symmetric functions}
$\sigma_l(T)$  of the $m$ variables in $T_j$  are defined by
$\sigma_0(T_j)=0$ and for
$1\le l\le m$, by
$$
\sigma_l(T_j)=\sum_{0\le r_1< r_2<\ldots<r_l\le
m-1}t_{j-r_1}t_{j-r_2}\ldots t_{j-r_l}.
$$
Let $C_m^l=\frac{m!}{l!(m-l)!}$ be the binomial coefficients, then the
monomials $e_l(x)=x^l$ can be written
$e_l=\sum_{i\in J}\theta_i^{(l)} B_i$, with
$\theta_i^{(l)}=\sigma_l(T_i)/C_m^l$, for $0\le l\le m$.
This is a direct consequence of Marsden's identity ([3], chapter IX).

*********************************

\sect{3. Differential QIs}

\noin
For all $j\in J$, we define $\psi_j(t)=\prod_{r\in \NN_m}(t_{j-r}-t)$
(thus $\psi_j\in \PP_{m}$ for all $j\in J$).
From de Boor and Fix [7] or de Boor ([5], chapter IX), we know that for any
$\tau\in \Sigma_j$,
the functionals
$$
\lambda_j(f)=\frac{1}{m!}\sum_{l=0}^{m} (-1)^{m-l}D^{m-l}\psi_j(\tau)D^lf(\tau)
$$
\noin
are dual functionals of B-splines, i.e. they satisfy, for all pairs
$(i,j)\in J\times J$
$$
\lambda_j(B_i)=\delta_{ij}.
$$
Therefore the {\sl differential quasi-interpolant} (abbr. DQI)
$$
Qf=\sum_{j\in J} \lambda_j(f) B_j
$$
satisfies $QB_j=B_j$ for all $j\in J$, i.e. $Q$ is a {\sl projector} on the
space $S_m(I,T)$.
In practice, it is interesting to choose
$\tau=\theta_j=\frac{1}{m}\sum_{s\in \NN_m} t_{j-s}=s_1(T_j)=\theta_j^{(1)}$.
However, the computation of $\lambda_j(f)$ needs the evaluation all
derivatives of polynomials $\psi_j$.
Another method consists in writing $Q$ in the form
$$
Qf=\sum_{i\in J} \tilde\lambda_i(f) B_i,
$$
whose coefficient functionals are defined by
$$
\tilde\lambda_i(f)=\sum_{l=0}^m a_l(\theta_i) \frac{D^lf(\theta_i)}{l!},
$$
and to impose that $Q$ be exact on monomials of degree at most $m$
$$
Qe_k=e_k \;\; for \;\; 0\le k\le m.
$$
Setting $\alpha_l(\theta_i)=\theta_i^{-l}a_l(\theta_i)$ and
$\beta_s(\theta_i)=\theta_i^{-s}\theta_i^{(s)}$, we obtain the following
system of linear
equations, for $0\le s\le m$:
$$
\sum_{l=0}^s C_s^l \alpha_l(\theta_i) =\beta_s(\theta_i).
$$
The solutions of this system are given by
$$
\alpha_s(\theta_i)=\sum_{l=0}^s (-1)^{s-l}C_s^l \beta_l(\theta_i).
$$
Thus we finally obtain
\proclaim Theorem 1.
The coefficients of the differential forms $\{\tilde\lambda_i(f), i\in J\}$,
are given, for $0\le s\le m$, by
$$
a_s(\theta_i)=\sum_{l=0}^s (-1)^{s-l}C_s^l \theta_i^{s-l}\theta_i^{(l)}.
$$

\noin
However, these DQIs need the values of derivatives of $f$, so they are not
very easy to use in
applications and we will not study them any more. Let us only give examples
of quadratic and cubic
DQIs.

\ms
\noin
{\bf Example 1:} {\sl Quadratic spline DQIs} (see also section 4 below).
In de Boor's form, we have for $\tau=\theta_j$:
$\lambda_j(f)=f(\theta_j)-\frac12 (\theta_j^2-\theta_j^{(2)}) D^2f(\theta_j)$,
where $\theta_j=\frac12(t_{j-1}+t_{j})$ and $\theta_j^{(2)}=t_{j-1}t_{j}$,
whence
$\theta_j^2-\theta_j^{(2)}=\frac14 (t_{j-1}-t_{j})^2$ and finally
$\lambda_j(f)=f(\theta_j)-\frac18 h_i^2 D^2f(\theta_j)$.
Theorem 1 gives
$a_0(\theta_j)=1,a_1(\theta_j)=0,a_2(\theta_j)=\theta_j^2-\theta_j^{(2)}$,
hence $\tilde\lambda_j(f)=\lambda_j(f)$.

\ms
\noin
{\bf Example 2:} {\sl Cubic spline DQIs}. In de Boor's form, we have for
$\tau=\theta_j$:
$\lambda_j(f)=f(\theta_j)-\frac12 (\theta_j^2-\theta_j^{(2)}) D^2f(\theta_j)
-\frac16 \psi_j(\theta_j)D^3f(\theta_j)$, with
$\theta_j=\frac13(t_{j-2}+t_{j-1}+t_{j})$
and $\theta_j^{(2)}=\frac13(t_{j-2}t_{j-1}+t_{j-1}t_{j}+t_{j-2}t_{j})$, whence
$\theta_j^2-\theta_j^{(2)}=\frac19 (h_{i-1}^2+h_{i-1}h_i+h_i^2)$ and
$\psi_j(\theta_j)=(t_j-\theta_j)(t_{j-1}-\theta_j)(t_{j-2}-\theta_j)$.
Finally, we obtain
$\lambda_j(f)=f(\theta_j)-\frac{1}{18}(h_{i-1}^2+h_{i-1}h_i+h_i^2)
D^2f(\theta_j)
-\frac{1}{162} (2h_{i-1}+h_i)(h_i-h_{i-1})(h_{i-1}+2h_i)D^3f(\theta_j)$.
Theorem 1 gives
$a_0(\theta_j)=1,a_1(\theta_j)=0,a_2(\theta_j)=\theta_j^2-\theta_j^{(2)},
a_3(\theta_j)=\theta_j^{(3)}-3\theta_j\theta^{(2)}+2\theta_j^3=
\frac{1}{27}(2h_{i-1}+h_i)(h_i-h_{i-1})(h_{i-1}+2h_i)$, whence
$\tilde\lambda_j(f)=\lambda_j(f)$.


\sect{4. Uniformly bounded discrete QIs exact on $\PP_2$}

\noin
It is now possible to derive {\sl discrete} QIs from the preceding DQIs by
replacing
the values of derivatives $D^lf(\theta_i)/l!$ of $f$ by divided differences
at the points
$\theta_r$ lying in
$\Sigma_i$. Doing this, we loose the property of projection on $S_m(I,T)$.
However, by choosing
conveniently the divided differences, we can obtain some families of dQIs
which are UB and
exact on specific subspaces of polynomials.

Let us construct for example a family of dQIs of degree $m$ which are {\sl
exact on} $\PP_2$.
We start from functionals which are truncations at order $2$ of those of
the preceding section:
$$
\lambda_j^{(2)}(f)=\frac{1}{m!}\sum_{l=0}^{2}
(-1)^{m-l}D^{m-l}\psi_j(\tau)D^lf(\tau).
$$
As $\psi_j(t)$ is of degree $m$, we obtain successively
$D^{m}\psi_j(\tau)=(-1)^{m} m!$,

\noindent
$D^{m-1}\psi_j(\tau)=(-1)^{m} m!(\tau-\theta_j)$,
$D^{m-2}\psi_j(\tau)=\frac12 (-1)^{m} m!(\tau^2-2\theta_j t+\theta_j^{(2)})$.

\noindent
More specifically, taking $\tau=\theta_j$, we get
$$
D^{m-1}\psi_j(\theta_j)=0,\;\; D^{m-2}\psi_j(\theta_j)=\frac12 (-1)^m
m!(\theta_j^{(2)}-\theta_j^2)
$$
and we obtain the DQI exact on $\PP_2$
$$
Q_2f=\sum_{j\in J} \lambda_j^{(2)}(f) B_j,
$$
whose coefficient functionals are given by
$$
\lambda_j^{(2)}(f)=f(\theta_j)-\frac12 (\theta_j^2-\theta_j^{(2)})
D^2f(\theta_j).
$$
We recall the expansion (se e.g.[9][17]):
$$
\bar\theta_j^{(2)}=\theta_j^2-\theta_j^{(2)}=\frac{1}{m^2
(m-1)}\sum_{(r,s)\in \NN_m^2,r<
s}(t_{j-r}-t_{j-s })^2>0.
$$
On the other hand, $\frac12 D^2f(\theta_j)$ coincide on the space $\PP_2$
with the
second order divided difference $[\theta_{j-1}, \theta_j,\theta_{j+1}]f$,
therefore the dQI defined by
$$
Q_2^*f=\sum_{j\in J} \mu_j^{(2)}(f) B_j,
$$
with coefficient functionals
$$
\mu_j^{(2)}(f)=f(\theta_j)-\bar\theta_j^{(2)} [\theta_{j-1}, \theta_j,
\theta_{j+1}]f,
$$
is also exact on $\PP_2$. Moreover, one can write
$$
\mu_i^{(2)}(f)=a_if_{i-1}+b_i f_i+c_i f_{i+1}
$$
with
$
a_i=-\bar\theta_i^{(2)}/\Delta \theta_{i-1}(\Delta \theta_{i-1}+\Delta
\theta_{i}),
\;\;
b_i=1+\bar\theta_i^{(2)}/\Delta \theta_{i-1}\Delta \theta_{i},
\quad
$

\noin
$c_i=-\bar\theta_i^{(2)}/\Delta \theta_{i}(\Delta \theta_{i-1}+\Delta
\theta_{i}).$
So, according to the introduction
$$
\Vert Q_2^* \Vert_{\infty}\le \max_{i\in J} (\vert a_i \vert+\vert b_i
\vert +\vert c_i \vert)
\le 1+2\max_{i\in J} \bar\theta_i^{(2)}/\Delta \theta_{i-1}\Delta \theta_{i}.
$$
The following theorem extends a result given for quadratic splines in
[13][22][23].
\ms
\noin
\proclaim Theorem 2.
For any degree $m$, the dQIs $Q_2^*$ are uniformly bounded. More
specifically, for all
partitions of $I$:
$$
\Vert Q_2^* \Vert_{\infty}\le [\frac12 (m+4)]
$$

\pf
We only give the proof for $m=2k+1$, the case $m=2k$ being similar. For the
sake of simplicity, we
take $j=k$, i.e. we shall determine an upper bound of the ratio
$$
N_k/D_k=\bar\theta_k^{(2)}/\Delta \theta_{k-1}\Delta \theta_{k}
$$
with
$$
N_k=\bar\theta_k^{(2)}=\frac{1}{m^2 (m-1)}\sum_{1\le r<s\le m} (t_r-t_s)^2.
$$
Setting $H=\sum_{i=1}^{m-1}h_i$, then we get a lower bound for the denominator
$$
D_k=\frac{1}{m^2}(t_{m-1}-t_0)(t_m-t_1)=\frac{1}{m^2}(h_0+H)(H+h_{m-1})\ge
\frac{H^2}{m^2}
$$
The numerator $N_k$ is composed of $k-1$ pairs of sums $(S_p,S'_p)$
$$
S_p=\sum_{s-r=p} (t_r-t_s)^2,\;\; S'_p=\sum_{s-r=k+p-1} (t_r-t_s)^2,
$$
for  $1\le p\le k-1$.
Both sums contain at most $p$ times the terms $h_i^2$ and $2h_i h_j$ ($i\ne
j$), hence we can
write
$S_p+S'_p\le 2p H^2$, which implies
$$
N_k\le \frac{2H^2}{(m-1)^2(m-2)}(1+2+\ldots k-1)=\frac{kS^2}{2(m-1)^2},
$$
so, we get
$$
N_k/D_k\le k/2,
$$
and finally, for $m=2k+1$ odd
$$
\Vert Q_2^* \Vert_{\infty}\le k+2=\frac12 (m+3)=[\frac12 (m+4)]
$$
For $m=2k$, we obtain respectively $D_k\ge\frac{S^2}{4k^2}$ and $N_k\le
\frac{S^2}{4(2k-1)}$, whence $N_k/D_k\le \frac{k^2}{2k-1}$, and finally
for $m=2k$ even
$$
\Vert Q_2^* \Vert\le k+2=\frac12 (m+4)=[\frac12 (m+4)]
$$
\eop

****************************************************************************
************

\sect{6. Existence and characterization of near-best discrete QIs}

\subsect{6.1.  Existence of near-best dQIs}
\noin
We  consider the following family of dQIs defined, for the sake of
simplicity, on $I=\RR$ endowed with an
arbitrary non-uniform increasing sequence of knots $T=\{t_i; i\in \ZZ\}$,
$$
Qf=Q_{p,q}f=\sum_{i\in\ZZ} \mu_i(f) B_i.
$$
Their coefficient functionals depend on $2p+1$ parameters, with $p\ge m$
$$
\mu_i(f)=\sum_{s=-p}^{p}\lambda_i(s)f(\theta_{i+s}),
$$
and they are exact on the space $\PP_q$, where $q\le \min(m,2p)$.
The latter condition is equivalent to $Q e_r=e_r$ for all monomials of
degrees $0\le r\le q$.
It implies that for all indices $i$, the parameters $\lambda_i(s)$ satisfy
the system of $q+1$
linear equations:
$$
\sum_{s=-p}^{p}\lambda_i(s)\theta_{i+s}^r=\theta_i^{(r)},\quad 0\le r\le q.
$$
The matrix $V_i\in \RR^{(q+1)\times (2p+1)}$ of this system, with coefficients
$V_i(r,s)=\theta_{i+s}^r$, is a Vandermonde
matrix of maximal rank $q+1$, therefore there are $2p-q$ {\sl free parameters}.
Denoting $b_i\in \RR^{q+1}$ the vector in the right hand side, with components
$b_i(r)=\theta_i^{(r)},\quad 0\le r\le q$, we consider the sequence of minimization problems,
for $i\in \ZZ$:
$$
\min \Vert \lambda_i \Vert_1,\quad V_i \lambda_i=b_i.
$$
We have seen in the introduction that $\nu_1^*(Q)=\max_{i\in\ZZ}\min\Vert
\lambda_i \Vert_1$ is an upper bound
of $\Vert Q_q \Vert_{\infty}$ which is easier to evaluate than the true
norm of the dQI.

\proclaim Theorem 3. The above minimization problems have always solutions,
which, in general,
are non unique.

\pf
The {\sl objective function being convex} and {\sl the domains being affine
subspaces},
these classical optimization problems have always solutions, in general non
unique.

\ms
\noin
We postpone to sections 7 and 8 the computation of some optimal solutions
in the case $q=2$.

\subsect{6.2. Characterization of optimal solutions}

\noin
For $b\in \RR^m$ and $A\in\RR^{m\times n}$, let us consider the
$l_1$-minimization problem
$$
(1) \quad \min \Vert r(a)\Vert_1, \quad r(a)=b-Aa.
$$
We recall the characterization of optimal solutions for $l_1$-problems
given in [30], chapter 6.
Define the sets
$$
Z(a)=\{1\le i\le m\vert \; r_i(a)=0\}
$$
$$
V(a)=\{v\in \RR^m ; \Vert v\Vert_{\infty}\le 1, \; v_i=sgn(r_i(a)) \; for
\; i\notin Z(a) \}
$$
\proclaim Theorem 4. $a^*$ is a solution of (1) if and only if there exists
a vector $v^*\in V(a^*)$
satisfying $A^T v^*=0$.

*****************************************************************

\sect{7. A general family of spline discrete QIs exact on $\PP_2$}

\noin
In this section, we restrict our study to the subfamily of spline dQIs
which are {\sl exact on $\PP_2$},
i.e; we choose $q=2$.
We shall try to characterize optimal solutions in the sense of theorem 3
with the help of theorem 4. Let
$$
Q_{p,2} f=\sum_{i\in \ZZ} \mu_i(f) B_i,
$$
where the coefficient functionals depend on $2p+1$ parameters
$$
\mu_i(f)=\sum_{r=-p}^{p}\lambda_i(r) f(\theta_{i+r}).
$$
We shall need the following sets of indices
$$
\bar K=\{-p,\ldots,p\},\quad K^*=\{-p,0,p\},\quad K=\bar K\backslash K^*.
$$
$$
K=K_1\cup K_2, \quad K_1=\{-p+1,\ldots,-1\},\quad K_2=\{1,\ldots,p-1\}.
$$
The three equations expressing the exactness of $Q_{p,2}$ on $\PP_2$ can be
written
$$
\lambda_i(-p)+\lambda_i(0)+\lambda_i(p)=1-\sum_{r\in K}\lambda_i(r)
$$
$$
\theta_{i-p}\lambda_i(-p)+\theta_{i}\lambda_i(0)+\theta_{i+p}
\lambda_i(p)=\theta_{i}-\sum_{r\in K}\theta_{r}\lambda_i(r)
$$
$$
\theta_{i-p}^2\lambda_i(-p)+\theta_{i}^2\lambda_i(0)+\theta_{i+p}^2\lambda_i
(p)=\theta_i^{(2)}-\sum_{r\in K}\theta_{r}^2\lambda_i(r)
$$
This system has a positive Vandermonde determinant
$$
V_i=V(\theta_{i-p},\theta_{i},\theta_{i+p})=(\theta_{i}-\theta_{i-p})(\theta
_{i+p}-\theta_{i})(\theta_{i+p}-\theta_{i-p}).
$$
Let us denote by
$$
(\lambda_i^*(-p),\lambda_i^*(0),\lambda_i^*(p))
$$
the unique solution of the above system with the right-hand side obtained
by taking $\lambda_i(r)=0$ for
all $r\in K$. Using Cramer's rule and the determinants $W_i(s)$ obtained by
replacing the column of $\theta_{i+s}$ in $V_i$
by this rhs, we obtain
$$
\lambda_i^*(-p)=W_i(-p)/V_i,\quad \lambda_i^*(0)=W_i(0)/V_i,\quad
\lambda_i^*(p)=W_i(p)/V_i.
$$
Then we can express the general solution of the above system in the form
$$
\lambda_i(-p)=\lambda_i^*(-p)-\sum_{r\in
K_1}\alpha_r\lambda_i(r)+\sum_{s\in K_2}\alpha_s\lambda_i(s)
$$
$$
\lambda_i(0)=\lambda_i^*(0)-\sum_{r\in K_1}\beta_r\lambda_i(r)-\sum_{s\in
K_2}\beta_s\lambda_i(s)
$$
$$
\lambda_i(p)=\lambda_i^*(p)+\sum_{r\in K_1}\gamma_r\lambda_i(r)-\sum_{s\in
K_2}\gamma_s\lambda_i(s)
$$
The various coefficients are quotients of Vandermonde determinants
\ms
$\alpha_r=V(\theta_{r},\theta_{i},\theta_{i+p})/V_i,\quad
\alpha_s=V(\theta_{i},\theta_{s},\theta_{i+p})/V_i,$
\ms
$\beta_r=V(\theta_{i-p},\theta_{r},\theta_{i+p})/V_i, \quad
\beta_s=V(\theta_{i-p},\theta_{s},\theta_{i+p})/V_i,$
\ms
$\gamma_r=V(\theta_{i-p},\theta_{r},\theta_{i})/V_i, \quad
\gamma_s=V(\theta_{i-p},\theta_{i},\theta_{s})/V_i.$
\ms
\noin
We denote by $Q_{p,2}^*$  the spline dQI whose coefficient functionals are
$$
\mu_i^*(f)=\lambda_i^*(-p)f(\theta_{i-p})+\lambda_i^*(0)f(\theta_{i})+
\lambda_i^*(p)f(\theta_{i+p}).
$$
In that case, an upper bound of the norm is $\max_{i\in \ZZ} \nu_i^*$ where
$$
\nu_i^*=\vert \lambda_i^*(-p) \vert+\vert \lambda_i^*(0) \vert+\vert
\lambda_i^*(p) \vert
$$
\proclaim Theorem 5. For all $p\ge m=$ degree of the spline, the infinite
norms of the spline
dQIs $Q_{p,2}^*$ are uniformly bounded  by $\frac{m+1}{m-1}$. This bound is
independent of $p$ and
of the sequence of knots $T$.

\pf
We have to find a good upper bound of
$$
\nu_i^*=\vert \lambda_i^*(-p)\vert+\vert
\lambda_i^*(0)\vert+\vert\lambda_i^*(p) \vert
$$
where, expanding the determinants, we have
$$
\lambda_i^*(-p)=-\bar\theta_i^{(2)}/(\theta_{i+p}-
\theta_{i-p})((\theta_i-\theta_{i-p})
$$
$$
\lambda_i^*(0)=1+\bar\theta_i^{(2)}/(\theta_{i+p}-\theta_i)((\theta_i-\theta
_{i-p})
$$
$$
\lambda_i^*(p)=-\bar\theta_i^{(2)}/(\theta_{i+p}-\theta_{i-p})((\theta_{i+p}
-\theta_i).
$$
We recall that $\theta_i=\frac1m\sum_{r=0}^{m-1}t_{i-r}$ and
$$
\bar\theta_i^{(2)}=\frac{1}{m^2(m-1)}\sum_{(r,s)\in\NN_m,r<s}(t_{i-r}-t_{i-s
})^2=\frac{S_1}{m^2(m-1)}
$$
We first compute
$$
\theta_i-\theta_{i-p}=\frac1m\sum_{r\in\NN_m}(t_{i-r}-t_{i-p-r})=\frac1m\sum
_{r\in\NN_m} \sum_{k=0}^p h_{i-k-r+1}=S_2/m.
$$
$$
\theta_{i+p}-\theta_i=\frac1m\sum_{r\in\NN_m}(t_{i+p-r}-t_{i-r})=\frac1m\sum
_{r\in\NN_m} \sum_{k=0}^p h_{i+k-r+1}=S_3/m,
$$
The proof being essentially the same for all $p\ge m$ and for all $i\in \ZZ$,
we can restrict our study to the cases $p=m$ and $i=m-1$. In that case, we get
$$
S_2=mh_{i-m+1}+\sum_{k=1}^{m-1} k(h_{}+h_{})\ge
S'_2=h_1+2h_2+\ldots+(m-1)h_{m-1},
$$
$$
S_3=mh_{i-m+1}+\sum_{k=1}^{m-1} k(h_{}+h_{})\ge
S'_3=(m-1)h_1+(m-2)h_2+\ldots+2h_{m-2}+h_{m-1}.
$$
Denoting, for $1\le k\le m-1$
$$
s_k=h_1+\ldots+h_k,
$$
and $s=s_{m-1}$, we get
$$
S'_2=s_1+s_2+\ldots+s_{m-1},\;\; S'_3=s+(s-s_1)+(s-s_2)+\ldots
(s-s_{m-2})=ms-S'_2
$$
whence
$$
S_2 S_3\ge S'_2S'_3=ms(s_1+s_2+\ldots+s_{m-1})-(s_1+s_2+\ldots+s_{m-1})^2.
$$
Now, we come back to $S_1$ and we shall prove that $S_1\le S'_2S'_3\le S_2S_3$.
$S_1$ can be written under the form
$$
S_1=\sum_{j=1}^{m-1}\sum_{}^{}h_{i-r+j}=\sum_{i=1}^{m-1}s_i^2+\sum_{j=1}^{m-
1}\sum_{i=j+1}^{m}(s_i-s_j)^2,
$$
from which we deduce
$$
S_1=(m-1)\sum_{i=1}^{m-1}s_i^2-2\sum_{j=1}^{m-1}s_j\sum_{i=j+1}^{m}s_i.
$$
Moreover, for all $1\le i\le m-1$, we have
$$
(m-1)s_i^2=ms_i^2-s_i^2\le ms_i s_{m-1}-s_i^2
$$
therefore, we obtain the result
$$
S_1\le (m-1)\sum_{i=1}^{m-1}s_i^2\le
ms_{m-1}\sum_{i=1}^{m-1}s_i-\sum_{i=1}^{m-1}s_i^2=S'_2S'_3\le S_2S_3.
$$
Finally, for all $i\in \ZZ$, we have
$$
\nu_i^*=1+\frac{2}{m-1}\frac{S_1}{S_2S_3}\le 1+\frac{2}{m-1}=\frac{m+1}{m-1},
$$
whence $\Vert Q_{p,2}^* \Vert_{\infty}\le \max_{i\in\ZZ}\nu_i^*\le
\frac{m+1}{m-1}.$
\eop

\noin
In the next section, we prove that  the QIs $Q_{p,2}^*$ are near-best in
the sense of section 6.

*********************************

\sect{8. A family of near-best spline discrete QIs}

For dQIs $Q_{p,2}$ depending on $p\ge m$ parameters, the coefficients (see
proof of theorem 5) are given by
$$
\lambda_i^*(-p)=-\bar\theta_i^{(2)}/(\theta_{i+p}-\theta_{i-p})((\theta_i-
\theta_{i-p})
$$
$$
\lambda_i^*(0)=1+\bar\theta_i^{(2)}/(\theta_{i+p}-\theta_i)((\theta_i-\theta
_{i-p})
$$
$$
\lambda_i^*(p)=-\bar\theta_i^{(2)}/(\theta_{i+p}-\theta_{i-p})((\theta_{i+p}
-\theta_i).
$$
\ms
\noin
Now, let us write the minimization problem of section 7 in Watson's form.
Denote
$$
\tilde\lambda_i=(\lambda_i(-p+1),\ldots,\lambda_i(-1),\lambda_i(1),\ldots,\l
ambda_i(p-1))^T\in\RR^{2p-2}
$$
$$
\lambda_i^*=(\lambda_i^*(-p),0,\ldots,\lambda_i^*(0),0,\ldots,\lambda_i^*(p)
)^T\in\RR^{2p+1}
$$
Let $A_i\in\RR^{(2p+1)\times(2p-1)}$ be the matrix with the following
coefficients (notations of section 7)
$$
For \; r\in K_1: A_i(-p,r)=\alpha_r, \; A_i(0,r)=\beta_r, \;
A_i(p,r)=-\gamma_r,
$$
$$
For \; s\in K_2:  A_i(-p,s)=-\alpha_s, \; A_i(0,s)=\beta_s, \;
A_i(p,s)=\gamma_s,
$$
$$
For \; r\in K_1: \; A_i(r,r')=0, \; r'\neq r,\; A_i(r,r)=-1, \quad
A_i(r,s)=0,\; s\in K_2,
$$
$$
For \; s\in K_2: \; A_i(s,r)=0, \; r\in K_1,\quad  A_i(s,s)=-1,\; A_i(s,s')=0,\; s'\neq s.
$$
Then, using these notations, we can write
$$
\Vert \lambda_i \Vert_1=\Vert \lambda_i^*-A_i \tilde\lambda_i \Vert_1
$$
\proclaim Theorem 6. Assume that the sequence of knots $T$ satisfies, for
all $i\in \ZZ$, the
following properties
$$
\theta_{i-1}+\theta_{i}\le \theta_{i-p}+\theta_{i+p}\le \theta_{i}+\theta_{i+1},
$$
then, for all $i\in \ZZ$, $\lambda_i^*$ is an optimal solution of the local
minimization problem
$\min \Vert \lambda_i \Vert_1$. Thus, for all $p\ge m$, the  spline dQIs
$Q_{p,2}^*$ are  near-best and  their infinite norms are uniformly bounded
by $\frac{m+1}{m-1}$.
This bound is independent of $p$ and of the sequence of
knots $T$.

\pf
According to Watson's theorem, we must find a vector $v^*\in\RR^{2p+1}$
satisfying
$$
\Vert v^* \Vert_{\infty}\le 1, \quad A_i^T v^*=0, \quad
v^*(r)=sgn(\lambda_i^*(r)) \;\; for \;\; r=-p,0,p.
$$

Let us choose
$$
v^*(-p)=-1,v^*(0)=1,v^*(p)=-1,
$$
$$
v^*(r)=-\alpha_r+\beta_r+\gamma_r,\;\; for \;\; r\in K_1,
$$
$$
v^*(s)=-\alpha_s+\beta_s+\gamma_s,\;\; for \;\; s\in K_2.
$$
Then it is easy to verify that the equations $A_i^T v^*=0$ are satisfied.
Moreover, the above expressions of $\lambda_i^*(r)$ for $r=-p,0,p$ with
$\bar\theta_i^{(2)}>0$
imply that $sgn(v^*(r))=sgn(\lambda_i^*(r))$ for $r=-p,0,p$.
It only remains to prove that, for $(r,s)\in K_1\times K_2$
$$
\vert v^*(r)\vert=\vert -\alpha_r+\beta_r+\gamma_r \vert\le 1,\quad
\vert v^*(s)\vert=\vert \alpha_s+\beta_s-\gamma_s \vert\le 1.
$$
As $\beta_r=1-\alpha_r+\gamma_r$ for $r\in K_1$ and
$\beta_s=1+\alpha_s-\gamma_s$ for $s\in K_2$,
it is equivalent to prove
$$
0\le \alpha_r-\gamma_r \le 1,\quad 0\le \gamma_s-\alpha_s \le 1, \;\; for
\;\; (r,s)\in K_1\times K_2
$$
We only detail the proof for $r\in K_1$, that for $s\in K_2$ being quite
similar.
Using the Vandermonde determinants, we get
$$
\alpha_r-\gamma_r=V_i^{-1}(\theta_{i}-\theta_{r})(\theta_{i+p}-\theta_{i-p})
[(\theta_{i+p}+\theta_{i-p})-(\theta_{r}+\theta_{i})],
$$
As $\theta_{i}-\theta_{r}\ge 0$ and $\theta_{i+p}-\theta_{i-p}\ge 0$, we
shall have $\alpha_r-\gamma_r\ge 0$
if and only if
$$
\theta_{r}+\theta_{i}\le \theta_{i+p}+\theta_{i-p}
$$
for all $r\in K_1$. However, since we have $\theta_{r}+\theta_{i}\le
\theta_{i-1}+\theta_{i}$,
there only remains the unique condition
$$
\theta_{i-1}+\theta_{i}\le \theta_{i-p}+\theta_{i+p}.
$$
The other inequality $\alpha_r-\gamma_r \le 1$ can be written
$$
(\theta_{i}-\theta_{r})[(\theta_{i+p}+\theta_{i-p})-(\theta_{r}+\theta_{i})]
\le (\theta_{i}-\theta_{i-p})(\theta_{i+p}-\theta_{i})
$$
Setting $\delta_1=\theta_{r}-\theta_{i-p}$,
$\delta_2=\theta_{i}-\theta_{r}$, and $\delta_3=\theta_{i+p}-\theta_{i}$,
the latter inequality can be written
$$
\delta_2(\delta_3-\delta_1)\le \delta_3 (\delta_2+\delta_1), \;\; or \;\;
\delta_1(\delta_2+\delta_3)\ge 0
$$
which is obviously satisfied.
For $s\in K_2$, the inequalities $0\le \gamma_s-\alpha_s \le 1$ are
satisfied if and only if
$$
\theta_{i-p}+\theta_{i+p}\le \theta_{i}+\theta_{i+1},
$$
whence the conditions on the sequence of knots.
\eop

{\bf Remark.} Theorem 6 imposes some conditions on the sequence of knots.
For quadratic splines,
we have studied arithmetic and geometric sequences: in both cases, the
higher is $p$,
the stronger are the conditions and for $p\to +\infty$, $T$ is closer and
closer to a uniform sequence.

***************************************************

\sect{9. Some applications}

\subsect{9.1. Approximation of functions}

\noin
When a spline dQI $Q$ is uniformly bounded independently of the partition,
we can apply a classical
result in approximation theory (see [5], Th 22, and [11], chapters 2 and 5):
$$
\Vert Qf-f\Vert_{\infty}\le (1+\Vert Q \Vert_{\infty}) d_{\infty}(f, \cal{S})
$$
where $\cal{S}$ is the space of splines. In particular, when $Q$ is exact
on the space $\PP_m$,
then for
$f\in C^{m+1}(I)$, one has
$$
\Vert Qf-f\Vert_{\infty}\le Ch^{m+1}\Vert f^{m+1}\Vert_{\infty}
$$
for some constant $C$ which does not depend on the given partition. Therefore spline dQIs give the best
possible  approximation order. More detailed results on error bounds are
given in [1],[13] and
[23].

\subsect{9.2. Quadrature formulas}

\noin
Approximating $\int_{I}f$ by $\int_{I}Q_2^*(f)$, where $Q_2^*(f)$ is the
quadratic spline dQI of section $8$, gives
rise to an interesting quadrature formula
$$
\int_{I}Q_2^*(f)=f_0\int_{I} B_0+\sum_{i=1}^n \mu_i(f)\int_{I}
B_i+f_{n+1}\int_{I} B_{n+1}
$$
As it is well known, $\int_{I} B_0=\frac{h_1}{3}$, $\int_{I}
B_{n+1}=\frac{h_n}{3}$ and $\int_{I}
B_i=\frac{h_{i-1}+h_i+h_{i+1}}{3}$ for $1\le i\le n$. This formula is exact
on $\PP_2$ , but in the
case of a uniform partition,(see [21]), it is exact on $\PP_3$ and provides
an interesting
complementary formula to Simpson's rule in the sense that, in general,
errors for both formulas have
opposite signs. This will be detailed in another paper, together with
applications to integral
equations.

\subsect{9.3. Pseudo-spectral methods}

\noin
One can approximate the first derivatives of a given function $f$ at the
data sites
$$
\Theta_n=\{\theta_0=t_0,\;\; \theta_i=\frac12(t_{i-1}+t_i),\;\; for \;\;
1\le i\le n,\;\;
\theta_{n+1}=t_n\}.
$$
by the derivatives of the quadratic spline dQI of section $8.1$
$$
Q^*_2 f=f(t_0)B_0+\sum_{i=1}^n \mu_i(f) B_i+f(t_n) B_{n+1}.
$$
For interior points $\theta_i,\;\; 3\le i\le n-2$, we obtain the general formula
$$
(Q^*_2 f)'(\theta_i)=\mu_{i-1}(f) B'_{i-1}(\theta_i)+\mu_i(f)
B'_i(\theta_i)+\mu_{i+1}(f)
B'_{i+1}(\theta_i)
$$
which can also be written, by setting $f_j=f(\theta_j)$:
$$
(Q^*_2 f)'(\theta_i)=\frac{1}{h_i}\{ -\sigma_i a_{i-1}f_{i-2}+[-\sigma_i
b_{i-1}+(\sigma_i-\sigma'_{i+1}) a_i]f_{i-1}
$$
$$
[-\sigma_i c_{i-1}+(\sigma_i-\sigma'_{i+1})b_i+a_{i+1}]f_i
+[(\sigma_i-\sigma'_{i+1})c_i+b_{i+1}]f_{i+1}+\sigma'_{i+1}c_{i+1}f_{i+2}\}
$$
For the first indices $0\le i\le 2$, the coefficients are modified
according to the convention
$h_0=0$, which gives $\sigma_0=0, \sigma'_0=1, \sigma_1=1$ and $\sigma'_1=0$.
We thus obtain
$$
(Q^*_2 f)'(\theta_0)=\frac{2}{h_1}\{ (a_1-1)f_0+b_1f_1+c_1f_2 \}
$$
$$
(Q^*_2f)'(\theta_1)=\frac{1}{h_1}\{(\sigma_2a_1-1)f_0+[\sigma_2b_1+\sigma'_2
a_2]f_1+
$$
$$
[\sigma_2c_1+\sigma'_2b_2]f_2+\sigma'_2c_2f_3\}
$$
$$
(Q^*_2f)'(\theta_2)=\frac{1}{h_2}\{-\sigma_2a_1f_0+[-\sigma_2b_1+(\sigma_2-\
sigma'_3)a_2]f_1+
[-\sigma_2 c_1+(\sigma_2-\sigma'_3) b_2+\sigma'_3 a_3]f_2
$$
$$
+[(\sigma_2-\sigma'_3) c_2+\sigma'_3 b_3]f_3+\sigma'_3 c_3 f_4\}
$$
In the same way, for the last
indices $n-1\le i\le n+1$, the coefficients are modified according to the
convention $h_{n+1}=0$ and
we obtain similar formulas for $(Q^*_2f)'(\theta_{n-1}),
(Q^*_2f)'(\theta_n)$ and
$(Q^*_2f)'(\theta_{n+1})$. In the case of a unifom partition, the formulas
are given in [21].
Concerning error estimates, it is rather easy to verify that $(Q^*_2
f)'(\theta_i)-f'(\theta_i)=O(h^2)$ where
$h=\max_{1\le i\le n}h_i$. A more detailed study will be done elsewhere.
These results can be used in pseudo-spectral methods, as described for
example in [12] and [29].

************************************************************

\References

\ref
D. Barrera, M.J. Iban\~ez, P. Sablonni\`ere:
Near-best discrete quasi-interpolants on uniform and nonuniform partitions.
In {\sl Curve and
Surface Fitting}, Saint-Malo 2002, A. Cohen, J.L. Merrien and L.L.
Schumaker (eds), Nashboro
Press, Brentwood (2003), 31-40.

\ref
D. Barrera, M.J. Iban\~ez, P. Sablonni\`ere, D. Sbibih:
Near-minimally normed univariate spline quasi-interpolants on uniform
partitions.
Pr\'epublication IRMAR 04-12, Universit\'e de Rennes, March 2004.

\ref
D. Barrera, M.J. Iban\~ez, P. Sablonni\`ere, D. Sbibih:
Near-best quasi-interpolants associated with H-splines on a three-direction
mesh.
\break
Pr\'epublication IRMAR 04-14, Universit\'e de Rennes, March 2004.

\ref
B.D. Bojanov, H.A. Hakopian, A.A. Sahakian:
{\sl Spline functions and multivariate interpolation}, Kluwer, Dordrecht 1993.

\ref
C. de Boor:
{\sl A practical guide to splines}, Springer-Verlag, New-York 2001.
(revised edition).

\ref
C. de Boor:
Splines as linear combinations of B-splines, a survey.
In: {\sl Approximation Theory II}, G.G. Lorentz et al. (eds), 1--47,
Academic Press, New-York 1976.

\ref
C. de Boor, G. Fix: Spline approximation by quasi-interpolants.
\JAT \quad {\bf 8} (1973), 19-45.

\ref
C. de Boor, K. H\"ollig, S. Riemenschneider:
{\sl Box-splines}, Springer-Verlag, New-York 1993.

\ref
G. Chen, C.K. Chui, M.J. Lai:
 Construction of real-time spline quasi-interpolation schemes,
\ATA {\bf 4} (1988), 61-75.

\ref
C.K. Chui:
{\sl Multivariate splines}, CBMS-NSF Regional Conference Series in Applied
Mathematics, vol. 54,
SIAM, Philadelphia 1988.

\ref
R.A. DeVore, G.G. Lorentz:
{\sl Constructive approximation}, Springer-Verlag, Berlin 1993.

\ref
B. Fornberg:
{\sl A practical guide to pseudospectral methods},
Cambridge University Press 1996.

\ref
M.J. Iba\~nez-P\'erez:
Cuasi-interpolantes spline discretos con norma casi minima : teoria y
aplicaciones.
Tesis doctoral, Universidad de Granada, 2003.

\ref
W.J. Kammerer, G.W. Reddien, R.S. Varga: Quadratic interpolatory splines.
{\sl Numer. Math.} {\bf 22}  (1974), 241--259

\ref
B.G. Lee, T. Lyche, L.L. Schumaker: Some examples of quasi-interpolants
constructed from local
spline projectors.  In {\sl Mathematical methods for curves and surfaces:
Oslo 2000},
T. Lyche and L.L. Schumaker (eds), Vanderbilt University Press, Nashville
(2001), 243-252.

\ref
T. Lyche, L.L. Schumaker:
Local spline approximation methods,
{\sl \JAT} {\bf 15} (1975), 294--325.

\ref
J.M. Marsden, I.J. Schoenberg:
An identity for spline functions with applications to variation diminishing
spline
approximation. {\sl \JAT} {\bf 3} (1970), 7--49.

\ref
J.M. Marsden:
 Operator norm bounds and error bounds for quadratic spline interpolation, In:
{\sl Approximation Theory}, Banach Center Publications, vol. {\bf 4}
(1979), 159--175.

\ref
G. N\"urnberger:
{\sl Approximation by spline function}, Springer-Verlag, Berlin 1989.

\ref
M.J.D. Powell: {\sl Approximation theory and methods}. Cambridge University
Press, 1981.

\ref
P. Sablonni\`ere: Bases de Bernstein et approximants splines. Th\`ese de
doctorat, Universit\'e
de Lille, 1982.

\ref
P. Sablonni\`ere: On some multivariate quadratic spline quasi-interpolants
on bounded domains.
In {\sl Modern developments in multivariate approximation}, W. Haussmann,
K. Jetter, M. Reimer,
J. St\" ockler (eds), ISNM Vol. 145, Birkh\"auser-Verlag, Basel (2003), 263-278.

\ref
P. Sablonni\`ere: Quadratic spline quasi-interpolants on bounded domains of
$\RR^d, d=1,2,3$.
{\sl Spline and radial functions},
Rend. Sem. Univ. Pol. Torino, Vol. {\bf 61} (2003), 61-78.

\ref
P. Sablonni\`ere: Recent progress on univariate and multivariate polynomial
or spline
quasi-interpolants.
Submitted to Proc. IBoMAT 2004,  Bommerholz, Germany (February 16-20, 2004).
Pr\'epublication IRMAR, Universit\'e de Rennes, March 2004.

\ref
P. Sablonni\`ere: Near-best univariate spline integral quasi-interpolants
on non-uniform partitions.
Pr\'epublication IRMAR, Universit\'e de Rennes, 2004 (in preparation).

\ref
I.J. Schoenberg:
{\sl Cardinal spline interpolation},
CBMS-NSF Regional Conference Series in Applied Mathematics, vol. 12,
SIAM, Philadelphia 1973.

\ref
I.J. Schoenberg:
{\sl Selected papers}, Volumes 1 and 2, edited by C. de Boor.
Birkh\"auser-Verlag, Boston 1988.

\ref
L.L. Schumaker:
{\sl Spline functions: basic theory}, John Wiley \& Sons, New-York 1981.

\ref
L. N. Trefethen:
{\sl Spectral methods in Matlab},
SIAM, Philadelphia, 2000.

\ref
G.G. Watson:
{\sl Approximation theory and numerical methods},
John Wiley and Sons, New-York, 1980.

******************************************************************

\Address

D. Barrera, M.J. Iba\~nez,
Departamento de Matem\'atica Aplicada,
Facultad de Ciencias, Universidad de Granada,
Campus de Fuentenueva,
18071 GRANADA, Spain.
{\tt dbarrera@ugr.es, mibanez@ugr.es}

\ms
 P. Sablonni\`ere, INSA de Rennes,
20 Avenue des Buttes de Co\"esmes,
 CS 14315, 35043 RENNES Cedex, France.
{\tt psablonn@insa-rennes.fr}

\ms
D. Sbibih,
D\'epartement de Math\'ematiques et Informatique,
Facult\'e des Sciences, Universit\'e Mohammed 1er,
40000 OUJDA, Marocco.
{\tt sbibih@sciences.univ-oujda.ac.ma}
\bye